\documentclass[12pt,twoside]{amsart}   
\usepackage{amsthm, amsmath,amssymb, amsfonts, fancyhdr}

\usepackage{graphicx,psfrag}

\usepackage{setspace}

\input xy
\xyoption{all}

\theoremstyle{plain}
\newtheorem{thm}{Theorem}[section]
\newtheorem{lem}[thm]{Lemma}
\newtheorem{prop}[thm]{Proposition}
\newtheorem{cor}[thm]{Corollary}

\theoremstyle{definition}
\newtheorem{defn}{Definition}[section]
\newtheorem{exmp}{Example}[section]

\newtheorem{notation}{Notation}[section]

\begin{document}

\title{A Poset Structure on Quasifibonacci Partitions}
\author{Hansheng Diao}

\maketitle

\pagestyle{myheadings} \markboth{\sc Hansheng Diao}{\sc
Quasifibonacci Partitions}

\begin{abstract}
In this paper, we study partitions of positive integers into
distinct quasifibonacci numbers. A digraph and poset structure is
constructed on the set of such partitions. Furthermore, we discuss
the symmetric and recursive relations between these posets. Finally,
we prove a strong generalization of Robbins' result on the
coefficients of a quasifibonacci power series.
\end{abstract}

\section{Introduction and Statement of Results}
\bigskip

Let $F_k$ denote the Fibonacci numbers (where we have shifted the
usual initial condition); i.e. $(F_1, F_2, \cdots) =
(1,2,3,5,8,\cdots)$.\\
Consider the formal power series
\begin{equation*}
\begin{split}
H(x) &=\prod_{k\geq1}(1-x^{F_k})\\
&=(1-x)(1-x^2)(1-x^3)(1-x^5)(1-x^8)\cdots\\
&=1-x-x^2+x^4+x^7-x^8+x^{11}-x^{12}-x^{13}+x^{14}+\cdots\\
\end{split}
\end{equation*}
Let $h_m$ be the coefficient of $x^m$ in $H(x)$. It is clear that
$h_m = h^+_m - h^-_m$, where $h^+_m$ is the number of partitions of
$m$ into an even number of distinct Fibonacci numbers, and $h^-_m$
is the number of partitions of $m$ into an odd number of distinct
Fibonacci numbers.

In [2], N. Robbins proved that $h_m \in \{-1, 0, 1\}$. In [1], F.
Ardila gave a simpler proof for Robbins' result by giving a
recursion on $h_m$.

In this paper, we consider \emph{quasifibonacci numbers}, which
serve as generalization of Fibonacci numbers.

\begin{defn}
Given a positive integer $N \geq 2$. A sequence $A_1, A_2, \ldots$
of positive integers is called \emph{quasifibonacci sequence of
level $N$} if
\begin{itemize}
\item $A_{k+N}= A_{k+N-1}+\cdots+A_k $ for all $k \in \mathbb{Z}^+$,
and
\item $A_k > A_{k-1}+\cdots+A_1$ for all $1\leq k \leq N$.
\end{itemize}
We also say that $A_1, A_2,\cdots$ are quasifibonacci numbers.
\end{defn}
In particular, (shifted) Fibonacci numbers and Lucas numbers are
quasifibonacci.

As we will see in this paper, Robbins' result can be generalized to
\emph{quasifibonacci numbers}. More precisely, we shall prove the
following theorem:
\begin{thm}
Let $A_1, A_2, \cdots$ be a quasifibonacci sequence of even level.
Consider the formal power series
\begin{equation*}
\begin{split}
H(x) &= \prod_{k\geq1}(1-x^{A_k})\\
&= (1-x^{A_1})(1-x^{A_2})(1-x^{A_3})(1-x^{A_4})\cdots\\
&=1+ \sum_{m\geq1}h_mx^m
\end{split}
\end{equation*}
Then $h_m \in \{-1, 0, 1\}$
\end{thm}

Similarly, we have $h_m = h^+_m - h^-_m$ where $h^+_m$ (resp.
$h^-_m$) is the number of partitions of $m$ into an
even (resp. odd) number of distinct quasifibonacci numbers.\\
\indent In this paper, we study the structure of the set of such
partitions. In fact, a digraph and poset structure on such sets will
be constructed in section 3. In section 4, we will unveil intrinsic
symmetry and recursive relations between these posets. Finally, as
an application, we shall prove Theorem 1.1 in section 5.

\section{Notations}
\bigskip
\begin{notation}The following notations will be used throughout the
paper.
\begin{itemize}
\item $\{0,1\}^{\omega}:= \{(a_1,a_2,\cdots) | a_i \in \{0,1\}, a_i=0 \textrm{ for all but finitely many }i\textrm{'s}\}$
\item Given a quasifibonacci sequence $A_1, A_2, \cdots$, define
\begin{equation*}
\begin{split}
S_n &= S_n(\{A_k\})\\
&:= \{(a_1,a_2,\cdots)\in \{0,1\}^{\omega} |
\sum_{i=1}^{\infty}a_iA_i = n, a_i\in \{0,1\} \}.\\
\end{split}
\end{equation*}
$S_n$ represents the set of partitions of $n$ into distinct
quasifibonacci numbers $A_k$.
\item For $a=(a_1,a_2,\cdots) \in \{0,1\}^{\omega}$, define \[s(a) := \sum_{i=1}^{\infty}a_iA_i.\]
We also say that $a$ is the \emph{representation} of $s(a)$.
\item For $k \geq N+1$, define
\begin{equation*}
\begin{split}
& A_{k,0} := A_k+\sum_{\substack{1\leq i < k-N-1 \\ N\nmid
i}}A_{k-N-1-i}\\
& A_{k,1} := A_k+\sum_{\substack{1\leq i < k-N \\ N\nmid
i}}A_{k-N-i}\\
& A_{k,2} := A_k+\sum_{\substack{1\leq i < k-N+1 \\ N\nmid
i}}A_{k-N+1-i}\\
\end{split}
\end{equation*}
\item For any $(a_1,a_2,\cdots) \in S_n$, define the \emph{length}
$l(a)$ to be the largest $i$ such that $a_i = 1$. (Abusing the
notation, we also identify $a$ with $(a_1,a_2,\cdots,a_{l(a)})$,
which is called the \emph{reduced representation}.)

\end{itemize}
\end{notation}

The next lemma gives some important arithmetic properties of the
$A_k$'s which will be used frequently throughout the paper.
\begin{lem}
Let $A_1, A_2, \cdots$ be a quasifibonacci sequence. Then
\begin{enumerate}
\item $A_{k+2}> A_k+A_{k-1}+\cdots+A_1$ for any $k \in \mathbb{Z}^+$
\item for $A_k \leq n < A_{k+1}$ and any $a \in S_n$, we have $l(a) \in \{k-1,k\}$
\item for $A_1+A_2+\cdots+A_{k-1} < n < A_{k+1}$ and any $a\in S_n$,
we have $l(a) = k$.
\item \[\sum_{\substack{1\leq i \leq k-1 \\ N\nmid i}} A_{k-i} < A_k\]
\item $A_k < A_{k,0} < A_{k,1} < A_{k,2} <A_{k+1}$
\item $A_1+A_2+\cdots+A_k > 2A_{k,0}$
\end{enumerate}

\end{lem}

\begin{proof}
\begin{enumerate}
\item By the definition of quasifibonacci numbers, it is clear that
$A_{k+2}> A_{k+1}\geq A_{k}+\cdots+A_1$ when $1\leq k\leq N$. The
case for $k>N$ follows immediately from induction and the following
inequality:
\[A_{k+2}=A_{k+1}+\cdots+A_{k-N+2}\geq A_{k+1}+A_k\]
\item It is clear that $l(a)\leq k$. Suppose $l(a)\leq k-2$. By the
previous lemma, $n \leq A_{k-2}+\cdots+A_1 < A_k$, a contradiction.
Therefore, $l(a)\geq k-1$.
\item This is straightforward.
\item Write $k-1 = Nk_1+r$ ($0\leq r\leq N-1$). Then
\begin{equation*}
\begin{split}
\sum_{\substack{1\leq i\leq k-1\\ N\nmid i}} A_{k-i} &=\sum_{j=1}^{k_1}\sum_{i=1}^{N-1}A_{k-jN+i} + \sum_{i=1}^r A_i\\
&=\sum_{j=1}^{k_1}(A_{k-jN+N}-A_{k-jN})+\sum_{i=1}^r A_i\\
&=A_k-A_{r+1}+\sum_{i=1}^r A_i\\
&< A_k\\
\end{split}
\end{equation*}
\item It is obvious that $A_k < A_{k,0} < A_{k,1} < A_{k,2}$.\\
By (4),
\begin{equation*}
\begin{split}
A_{k,2} &= A_k+\sum_{\substack{1\leq i < k-N+1 \\ N\nmid
i}}A_{k-N+1-i}\\
& < A_k+A_{k-N+1}\\
& \leq A_k +A_{k-1}\\
& \leq A_{k+1}\\
\end{split}
\end{equation*}
\item By (4), we have
\begin{equation*}
\begin{split}
2A_{k,0}&=A_k + A_k +\sum_{\substack{1\leq i < k-N-1 \\ N\nmid
i}}A_{k-N-1-i} +\sum_{\substack{1\leq i < k-N-1 \\ N\nmid
i}}A_{k-N-1-i}\\
& < A_k +(A_{k-1}+\cdots+A_{k-N}) +\sum_{\substack{1\leq i <
k-N-1 \\ N\nmid i}}A_{k-N-1-i} + A_{k-N-1}\\
&\leq A_1 + A_2 +\cdots +A_k\\
\end{split}
\end{equation*}
\end{enumerate}
\end{proof}

\section{A Digraph and Poset Structure on $S_n$}
\bigskip
For each $n\geq 1$, we construct a digraph $G_n :=G_n(\{A_k\})$ in
the following way:
\begin{enumerate}
\item Set $V(G_n) = S_n$. In particular, set $G_n = \emptyset$
if $S_n = \emptyset$.
\item For $a=(a_1,a_2,\cdots)$, $b=(b_1,b_2,\cdots) \in S_n$,
let $(a,b) \in E(G_n)$ if there exists $k \in \mathbb{Z}^+$ such
that
\begin{itemize}
\item $a_{k+N}=1$, $a_k = a_{k+1}=\cdots= a_{k+N-1}=0$
\item $b_{k+N}=0$, $b_k=b_{k+1}=\cdots= b_{k+N-1}=1$
\item $a_t = b_t$ for all $t\notin \{k, k+1,\cdots, k+N\}$
\end{itemize}
(Here $(u, v)$ represents the directed edge $u\rightarrow v$.)
\end{enumerate}
\bigskip

The digraph structure induces a natural partial order on $S_n$ as follows:\\
\begin{center}
\emph{For $a,b \in S_n$, set $a\geq b$ if there exists a path
in $G_n$ from $a$ to $b$. }
\end{center}
\bigskip
In other words, $a$ covers $b$ if and only if $(a, b) \in E(G_n)$.

This makes $S_n$ into a poset $P_n$. The following examples show the
corresponding Hasse diagrams for $P_n(\{F_k\})$ and $P_n(\{L_k\})$,
where $\{F_k\}$ and $\{L_k\}$ denote the Fibonacci numbers and Lucas
numbers, respectively.

\begin{exmp}
First 24 Hasse diagrams for $P_n(\{F_k\})$:
\begin{figure}[h]

  \begin{center}
  \psfrag{p1}{$P_1$}
  \psfrag{p2}{$P_2$}
  \psfrag{p3}{$P_3$}
  \psfrag{p4}{$P_4$}
  \psfrag{p5}{$P_5$}
  \psfrag{p6}{$P_6$}
  \psfrag{p7}{$P_7$}
  \psfrag{p8}{$P_8$}
  \psfrag{p9}{$P_9$}
  \psfrag{p10}{$P_{10}$}
  \psfrag{p11}{$P_{11}$}
  \psfrag{p12}{$P_{12}$}
  \psfrag{p13}{$P_{13}$}
  \psfrag{p14}{$P_{14}$}
  \psfrag{p15}{$P_{15}$}
  \psfrag{p16}{$P_{16}$}
  \psfrag{p17}{$P_{17}$}
  \psfrag{p18}{$P_{18}$}
  \psfrag{p19}{$P_{19}$}
  \psfrag{p20}{$P_{20}$}
  \psfrag{p21}{$P_{21}$}
  \psfrag{p22}{$P_{22}$}
  \psfrag{p23}{$P_{23}$}
  \psfrag{p24}{$P_{24}$}
  \includegraphics[width=12cm]{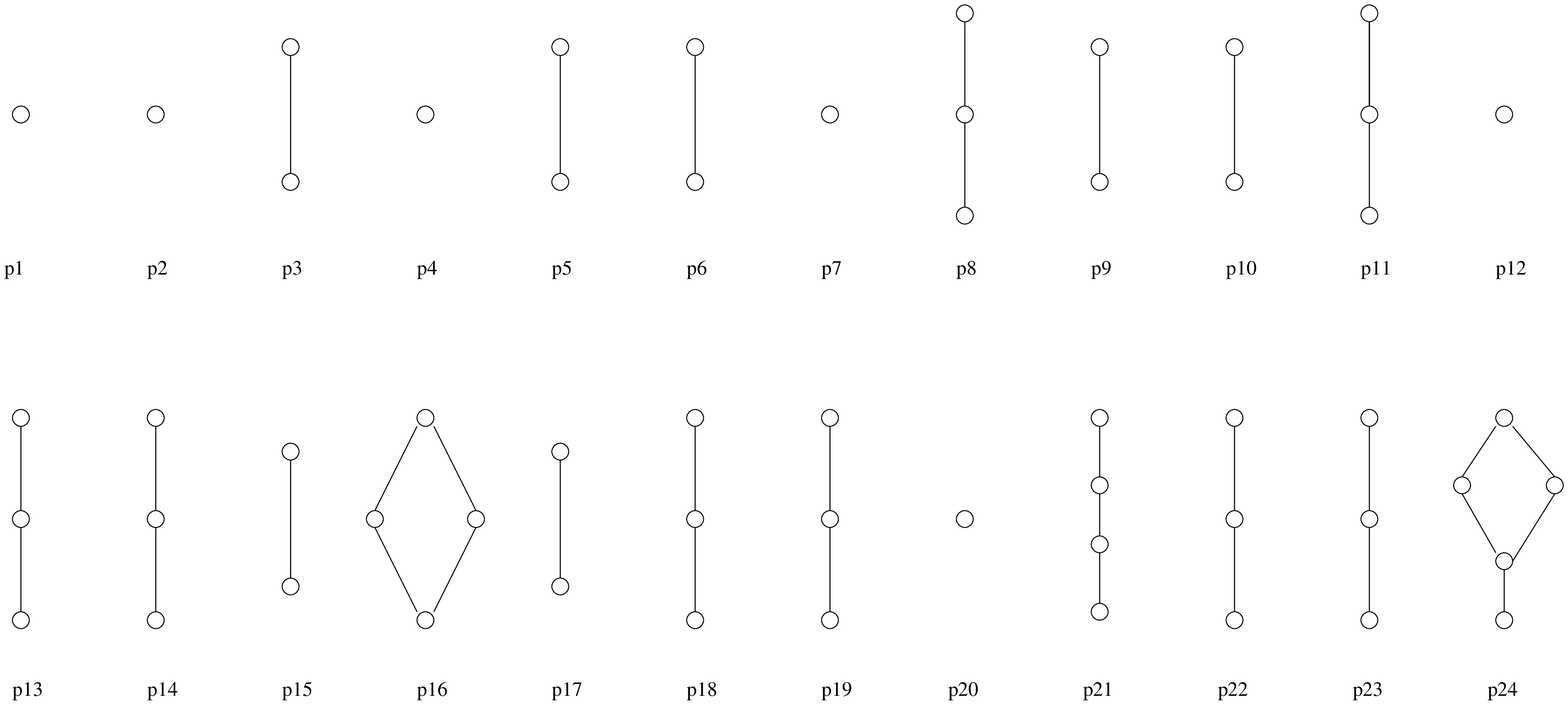}\\
  \end{center}
\end{figure}

\end{exmp}

\begin{exmp}
First 24 Hasse diagrams for $P_n(\{L_k\})$:
\begin{figure}[h]

  \begin{center}
  \psfrag{p1}{$P_1$}
  \psfrag{p2}{$P_2$}
  \psfrag{p3}{$P_3$}
  \psfrag{p4}{$P_4$}
  \psfrag{p5}{$P_5$}
  \psfrag{p6}{$P_6$}
  \psfrag{p7}{$P_7$}
  \psfrag{p8}{$P_8$}
  \psfrag{p9}{$P_9$}
  \psfrag{p10}{$P_{10}$}
  \psfrag{p11}{$P_{11}$}
  \psfrag{p12}{$P_{12}$}
  \psfrag{p13}{$P_{13}$}
  \psfrag{p14}{$P_{14}$}
  \psfrag{p15}{$P_{15}$}
  \psfrag{p16}{$P_{16}$}
  \psfrag{p17}{$P_{17}$}
  \psfrag{p18}{$P_{18}$}
  \psfrag{p19}{$P_{19}$}
  \psfrag{p20}{$P_{20}$}
  \psfrag{p21}{$P_{21}$}
  \psfrag{p22}{$P_{22}$}
  \psfrag{p23}{$P_{23}$}
  \psfrag{p24}{$P_{24}$}
  \psfrag{empty}{$\emptyset$}
  \includegraphics[width=12cm]{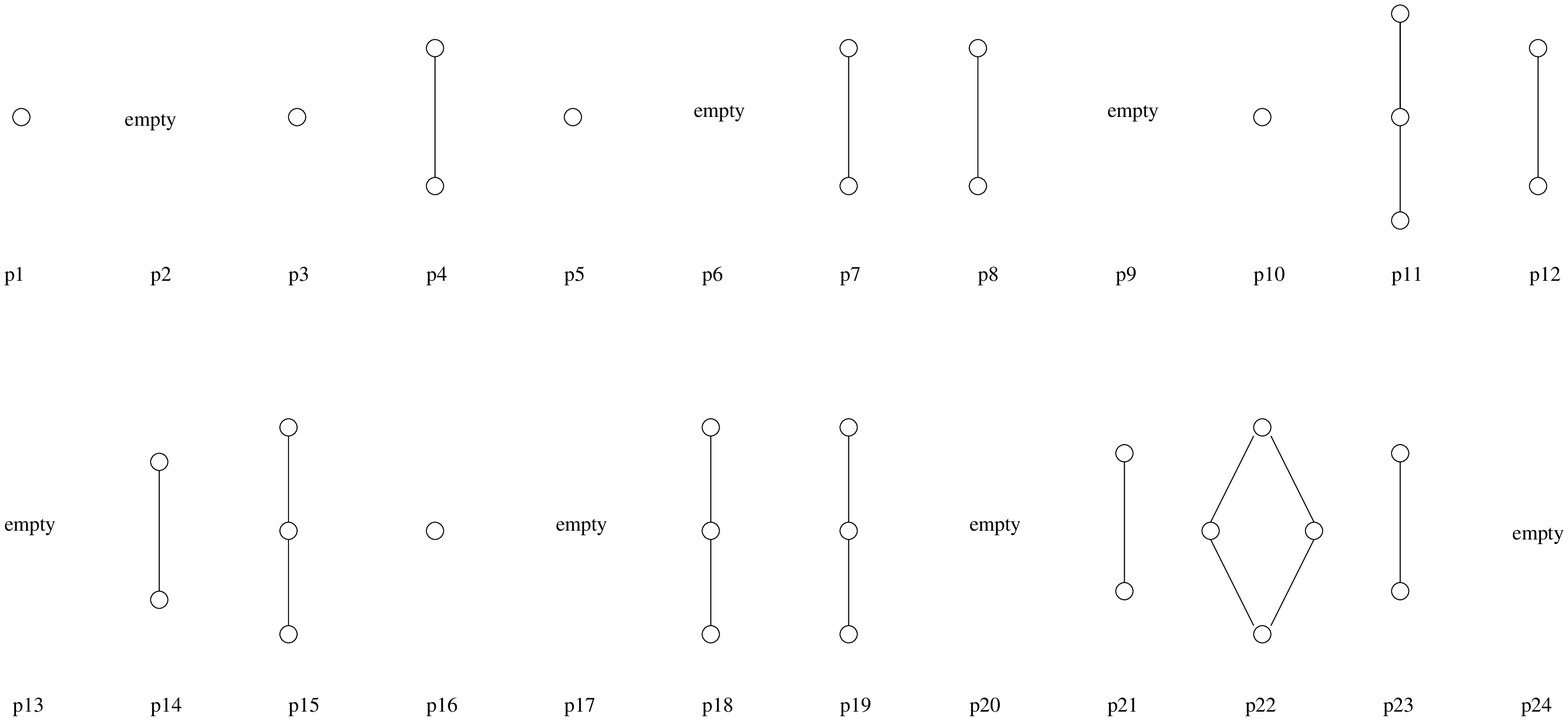}\\

  \end{center}
\end{figure}

\end{exmp}

\noindent For $A_k \leq n < A_{k+1}$, define
\begin{equation*}
\begin{split}
&T_n := \{a \in S_n | l(a) =k \}\\
&R_n := \{a \in S_n | l(a) = k-1\}\\
\end{split}
\end{equation*}
By Lemma 2.1(2), we have $S_n = T_n \cup R_n$. \\
Furthermore, let
$U_n$, $D_n$ denote the subposet of $P_n$ restricted on the vertex
set $T_n$ and $R_n$, respectively. (Abusing the notation, $U_n$,
$D_n$ also denote the corresponding subdigraphs of $G_n$.)

\begin{exmp}The figure below shows that $P_{24}$ can be decomposed into $U_{24}$ and $D_{24}$.
\begin{figure}[h]
\begin{center}
\psfrag{p24}{$P_{24}=P_{24}(\{F_k\})$} \psfrag{u24}{$U_{24}$}
\psfrag{d24}{$D_{24}$}
\includegraphics[width=13cm]{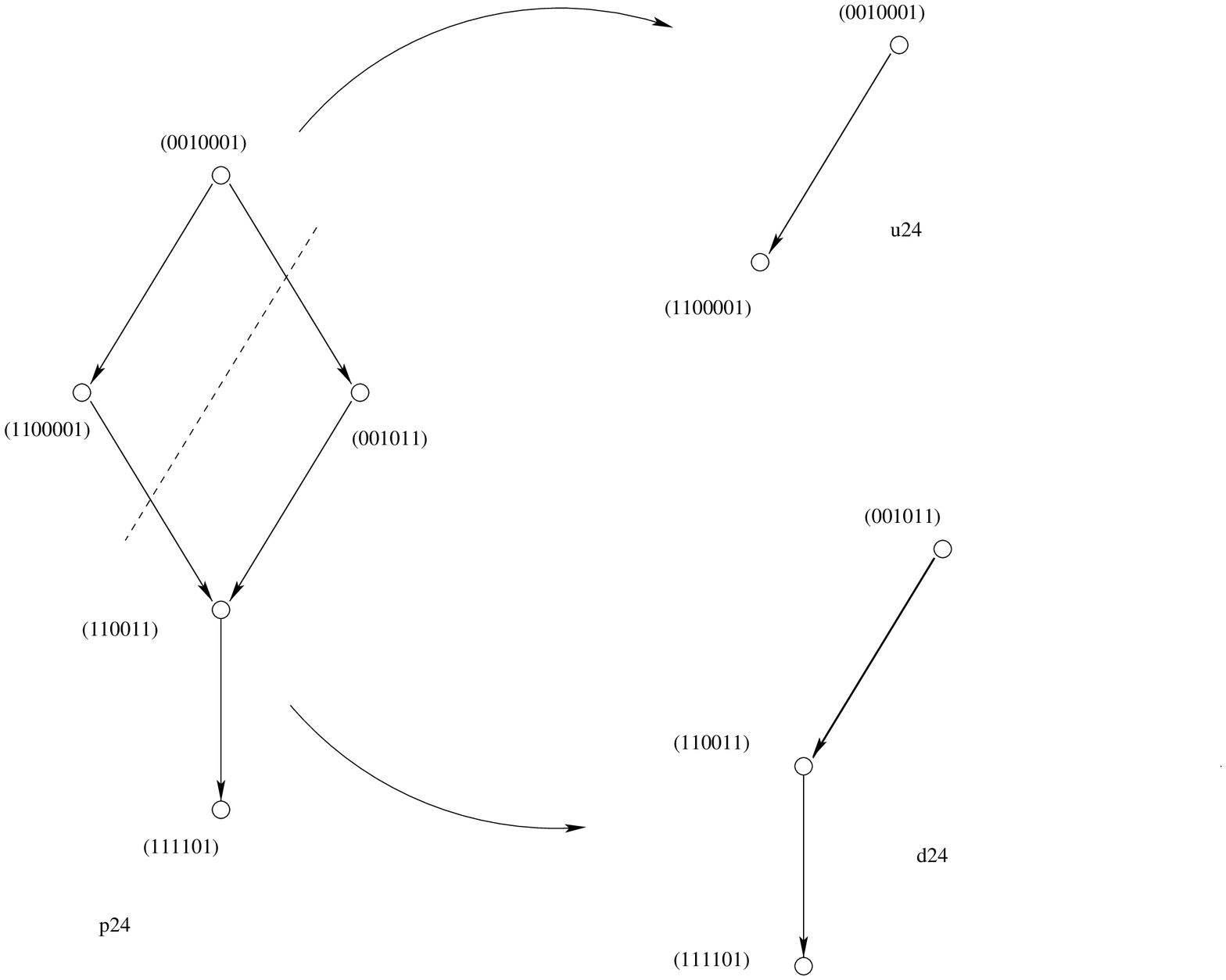}\\

\end{center}
\end{figure}

\end{exmp}

Now we study the structure of the posets $P_n$ in detail.
\begin{prop}
Suppose $S_n \neq \emptyset$. Then there is a unique maximal element
$\hat{1}$ in $P_n$. More precisely, $\hat{1}$ is the only element in
$S_n$ which does
not contain $N$ consecutive 1's.\\
Furthermore, if $A_k\leq n < A_{k+1}$, then $l(\hat{1}) = k$.

\end{prop}

\begin{proof}
Suppose $a\in S_n$ is a maximal element in $P_n$. It is clear that
$a$ does not contain $N$ consecutive 1's. (Otherwise, assume that
$a_k=a_{k+1}=\cdots=a_{k+N-1}=1$ and $a_{k+N}=0$. Set $b =
(b_1,b_2,\cdots)$ where $b_k=b_{k+1}=\cdots=b_{k+N-1}=0$,
$b_{k+N}=1$ and $b_t=a_t$ for $t\notin \{k,k+1,\cdots,k+N\}$. Then
$b \in P_n$ but $b>a$.)\\
Now we show the uniqueness. Suppose both $a$ and $a'$ are maximal
elements in $P_n$. Let $k$ be the largest index such that $a_k \neq
a'_k$. Without loss of generality, assume that $a_k=0$, $a'_k=1$ and
$a_t=a'_t=0$ for $t>k$. Since $a$ is maximal, it does not contain
$N$ consecutive 1's. It follows that
\[s(a) \leq \sum_{\substack{1\leq i\leq k-1\\ N\nmid i}} A_{k-i} \]
which is the largest possible value of
length $k-1$ with no $N$ consecutive 1's.\\
However, Lemma 2.1(4) gives \[s(a) \leq \sum_{\substack{1\leq i \leq
k-1 \\ N\nmid i}} A_{k-i} < A_k \leq s(a')\] a contradiction. Hence
the maximal element is unique. Denote it by
$\hat{1}$.\\
Now suppose $A_k \leq n < A_{k+1}$. By Lemma 2.1(2), $l(\hat{1}) \in
\{k, k-1\}$. As proved above, any element of length $k-1$ with no
$N$ consecutive 1's is smaller than $A_k$. So we must have
$l(\hat{1}) = k$.
\end{proof}

\begin{cor}
Suppose \[n> \sum_{\substack{1\leq i \leq k-1 \\ N\nmid i}}
A_{k-i}\] and $S_n \neq \emptyset$. Then $n \geq A_k$.
\end{cor}

\begin{proof}
Assume that \[\sum_{\substack{1\leq i \leq k-1 \\ N\nmid i}} A_{k-i}
< n < A_k.\] Then we must have $l(\hat{1}) \leq k-1$. Thus\[n =
s(\hat{1}) < \sum_{\substack{1\leq i \leq k-1 \\ N\nmid i}}
A_{k-i},\] a contradiction.
\end{proof}

Similarly, we have the following dual result:
\begin{prop}
Suppose $S_n \neq \emptyset$. Then there is a unique minimal element
$\hat{0}$ in $P_n$. More precisely, $\hat{0}$ is the only element in
$S_n$ which does not contain $N$ consecutive 0's in reduced
representation.

\end{prop}
However, we don't have $l(\hat{0})=k-1$ in general.

Actually, we will show that $P_n$ is a modular lattice in Section 5.

\section{Symmetry and Recursions}

If we viewing $P_1, P_2, \cdots$ as a sequence, then there exists
local symmetry relations between the posets. For instance, in
Example 3.1, the posets are central symmetric from $P_7$ to
$P_{12}$, and from $P_{12}$ to $P_{20}$. In general, similar
symmetry appears for all
quasifibonacci sequences.\\
In order to describe this special symmetry relation, we recall the
definition of \emph{dual posets}.
\begin{defn}
Two posets $P$, $Q$ are \emph{dual posets} to each other if there
exists an order-reversing bijection $\phi : P \rightarrow Q$ whose
inverse is also order-reversing; that is
\[x\leq y \textrm{ in } P \Leftrightarrow  \phi(y)\leq \phi(x)
\textrm{ in } Q\]
\end{defn}

\begin{prop}
For $A_k\leq n < A_{k+1}$, let $n' = A_1+A_2+\cdots+A_k-n$. Then
$P_n$ is dual to $P_{n'}$.
\end{prop}

\begin{proof}
Define $\phi: P_n \rightarrow P_{n'}$ by setting $(a_1,a_2,\cdots,
a_k) \mapsto (1-a_1,1-a_2,\cdots,1-a_k)$. (Note that $a_k$ is not
necessarily nonzero)\\
It is easy to check that $(a, b) \in E(G_n)$ if and only if
$(\phi(b),\phi(a))\in E(G_{n'})$. Hence $P_n$ is dual to $P_{n'}$
via $\phi$.

\end{proof}

Being more careful, we can derive similar symmetry on $U_n$ and
$D_n$.
\begin{prop}
For $A_k\leq n <A_k+A_{k-1}$, let $n'=A_1+A_2+\cdots+A_k-n$. Then
$U_n$ is dual to $D_{n'}$, and $D_n$ is dual to $U_{n'}$.
\end{prop}

\begin{proof}
It suffices to show that $\phi(U_n) = D_{n'}$ and $\phi(D_n) =
U_{n'}$.\\
Let $a = (a_1,a_2,\cdots,a_k) \in U_n$, ($a_k =1$). Since $n - A_k <
A_{k-1}$, $a_{k-1} = 0$. Thus $1-a_k = 0$, $1-a_{k-1}=1$.
Hence $l(\phi(a)) = k-1$.\\
On the other hand, for any $b = (b_1,b_2,\cdots,b_k) \in D_n$, we
must have $1-b_k = 1$. So $l(\phi(b)) = k$. \\
Therefore, $\phi(U_n) = D_{n'}$ and $\phi(D_n) = U_{n'}$, as
desired.
\end{proof}

Other than symmetry, there are intrinsic recursive relations in the
poset sequence $\{P_n\}$. In order to describe the recursion
clearly, we introduce the following notations.

\begin{notation}
Suppose $a=(a_1,a_2,\cdots)$, $b=(b_1,b_2,\cdots) \in
\{0,1\}^{\omega}$ satisfies $a_t = 0$ whenever $b_t=1$. Define
\[a+b := (a_1+b_1, a_2+b_2, \cdots)\] Similarly, suppose
$a=(a_1,a_2,\cdots)$, $b=(b_1,b_2,\cdots) \in \{0,1\}^{\omega}$
satisfies $a_t = 1$ whenever $b_t =1$. Define \[a-b := (a_1-b_1,
a_2-b_2,\cdots)\]
\end{notation}

\begin{notation}
Let $P$ be a finite subset of $\{0,1\}^{\omega}$. The poset
structure on $P$ is defined as usual. Let $b=(b_1,b_2,\cdots, b_k)
\in \{0,1\}^{\omega}$. Suppose that for any $a=(a_1,a_2,\cdots) \in
P$, we have $a_t=0$ whenever $b_t=1$. Define
\[P+b=\{a+b\, |\, a \in P\},\]
regarded as a poset, with the natural partial order.

\end{notation}

\begin{exmp}
The following figure gives an example for $P = P_{26}(\{F_k\})$ and
$b= (10000001)$:
\begin{figure}[h]
\begin{center}
\psfrag{P}{$P$} \psfrag{Pb}{$P$+(10000001)}
\includegraphics[width=12cm]{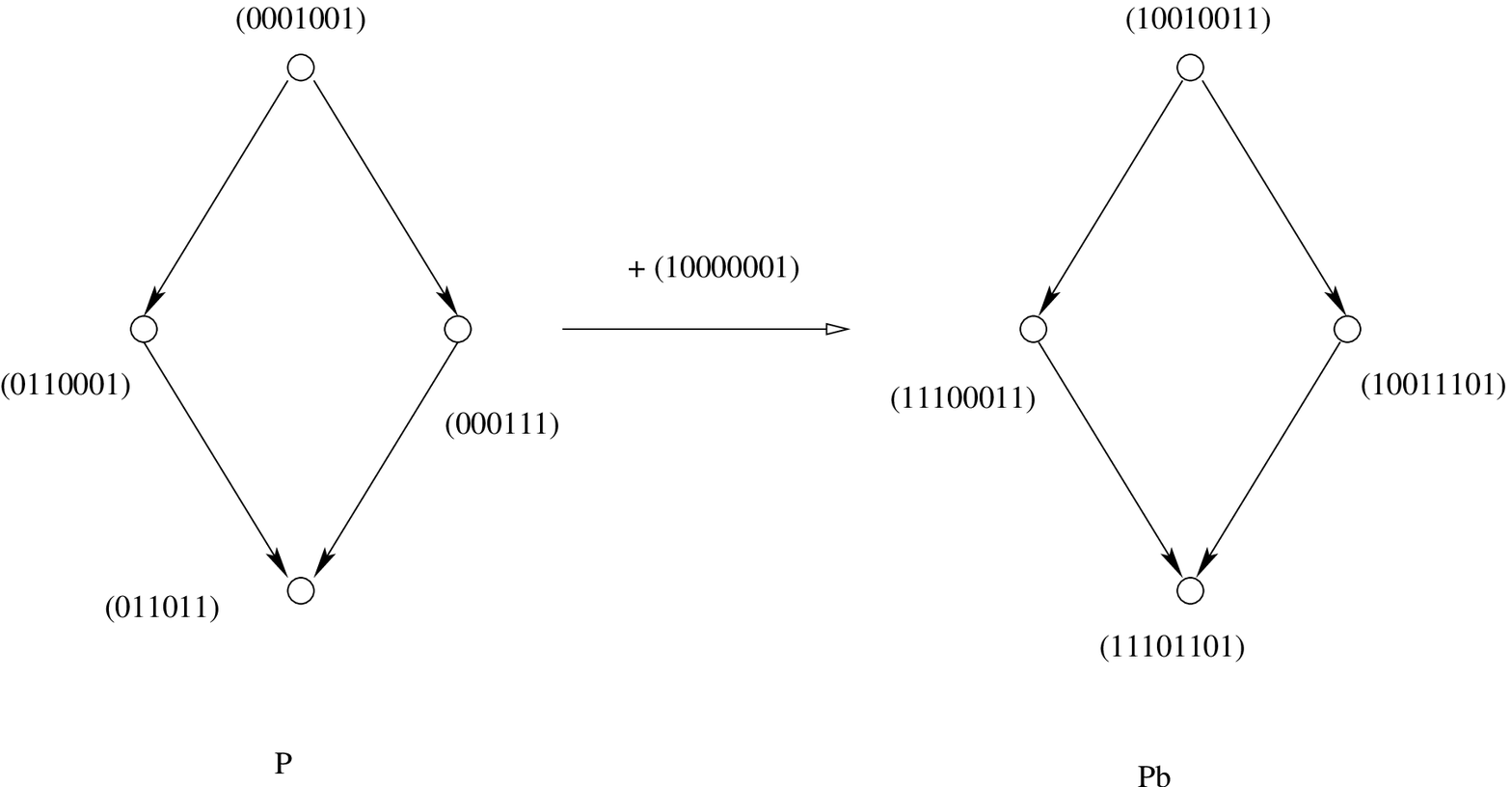}\\

\end{center}
\end{figure}

\end{exmp}

We can also define the addition between the posets:
\begin{notation}
Let $P$, $Q$ be disjoint finite subsets of $\{0,1\}^{\omega}$ with
natural digraph and poset structure. Let $P_1$, $Q_1$ be subposets
of $P$, $Q$, respectively. Suppose there is a bijection $\psi : P_1
\rightarrow Q_1$ such that $a$ covers $\psi(a)$ for all $a \in
P_1$.\\
Define
\[P\hat{+}Q = P \underset{(P_1,Q_1)}{\hat{+}}Q = P \cup Q \]
regarded as a digraph with
\[V(P\hat{+}Q) = V(P)\cup V(Q)\]
\[E(P\hat{+}Q) = E(P)\cup E(Q) \cup \{(a,\psi(a)) | a\in P_1\}\]
We also treat $P\hat{+}Q$ as a poset if the partial order determined
by the directed edges in $E(P\hat{+}Q)$ is exactly the natural one.
\end{notation}
In particular, $P_n = U_n\hat{+}D_n$.

\begin{figure}[h]
\begin{center}
\psfrag{p}{$P$} \psfrag{q}{$Q=Q_1$} \psfrag{pq}{$P\hat{+}Q$}
\psfrag{+}{$\hat{+}$} \psfrag{p1}{$P_1$} \psfrag{q1}{$Q_1$}
\includegraphics[width=13cm]{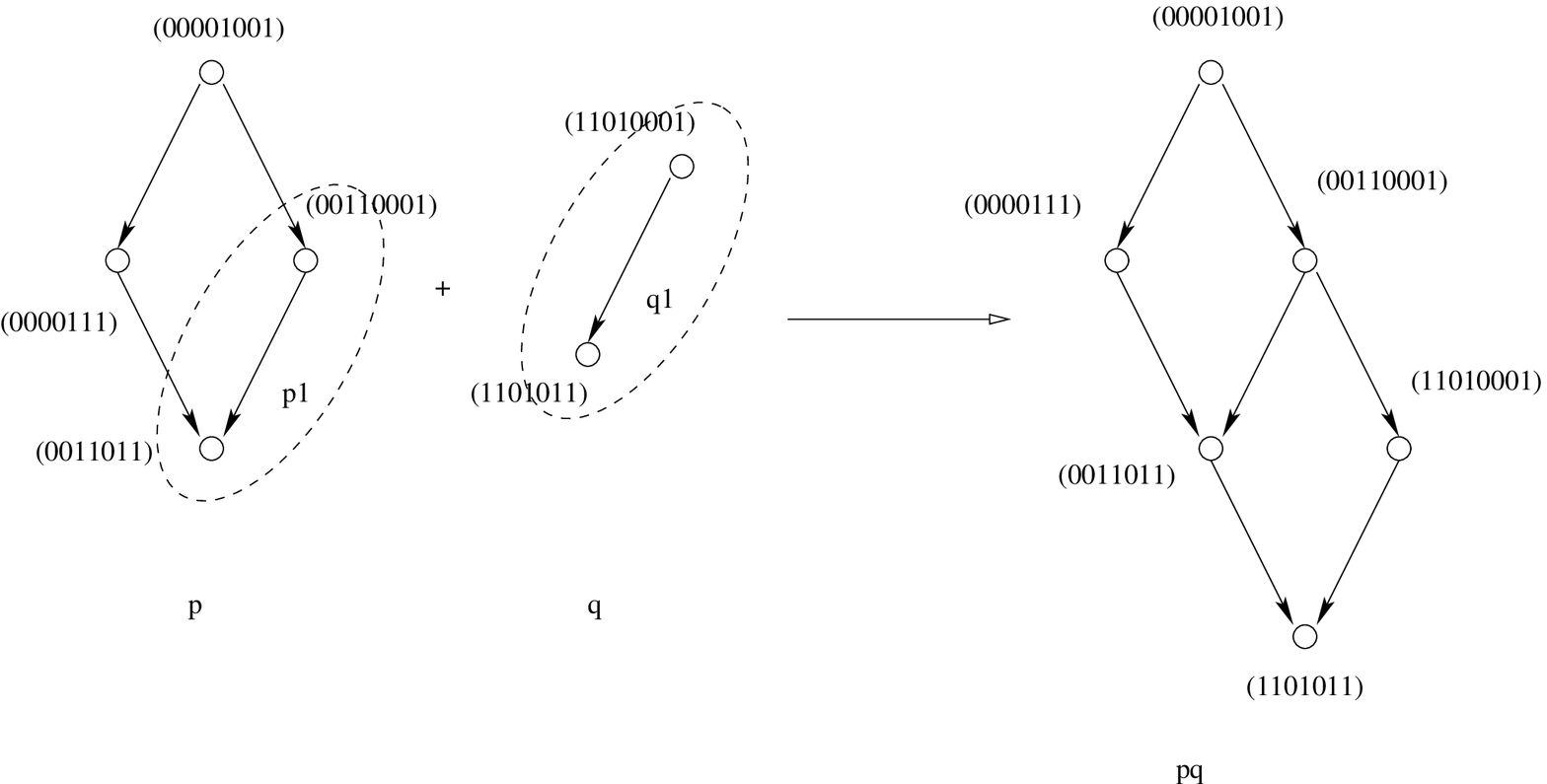}\\

\end{center}
\end{figure}

From now on, let $\tau_k$ denote the only element in
$\{0,1\}^{\omega}$ with the $k$-th entry being 1 and the others
being 0; i.e., $\tau_k = (00\cdots01) \in S_{A_k}$.\\
Let $\eta_k$ denote the only element $(a_1,a_2,\cdots, a_{k-1})$ in
$\{0,1\}^{\omega}$ with $a_{k-1}=a_{k-2}=\cdots=a_{k-N}=1$ and
$a_{k-N-1}=\cdots=a_1=0$; i.e., $\eta_k = (00\cdots0011\cdots11) \in
S_{A_k}$.

Now we describe the recursion explicitly. We show that each $P_n$
with $A_k\leq n < A_{k+1}$ can be expressed in terms of $P_1,
P_2,\cdots, P_{A_k-1}$.

\begin{prop}
If $A_{k,2}< n < A_{k+1}$, then $P_n$ is isomorphic to $P_{n-A_k}$.
More precisely, we have $P_n = P_{n-A_k}+\tau_k$, $U_n =
U_{n-A_k}+\tau_k$ and $D_n =D_{n-A_k}+\tau_k$.
\end{prop}

\begin{prop}
If $A_{k,1} < n \leq A_{k,2}$, then
\begin{equation*}
\begin{split}
P_n &= (P_{n-A_k}+\tau_k)\hat{+}(D_{n-A_k}+\eta_k)\\
&= (P_{n-A_k}+\tau_k)\underset{(D_{n-A_k}+\tau_k,
D_{n-A_k}+\eta_k)}{\hat{+}}(D_{n-A_k}+\eta_k)\\
\end{split}
\end{equation*}
$U_n = P_{n-A_k}+\tau_k$ and $D_n =D_{n-A_k}+\eta_k$.
\end{prop}

\begin{prop}
If $A_{k,0}< n \leq A_{k,1}$, then
\begin{equation*}
\begin{split}
P_n &= (P_{n-A_k}+\tau_k)\hat{+}(P_{n-A_k}+\eta_k)\\
&= (P_{n-A_k}+\tau_k)\underset{(P_{n-A_k}+\tau_k, P_{n-A_k}+\eta_k)}{\hat{+}}(P_{n-A_k}+\eta_k)\\
\end{split}
\end{equation*}
$U_n = P_{n-A_k}+\tau_k$ and $D_n=P_{n-A_k}+\eta_k$.
\end{prop}

The only remaining case is $A_k \leq n \leq A_{k,0}$. By Proposition
4.1, $P_n$ is dual to $P_{n'}$ where
$n'=A_1+A_2+\cdots+A_k-n$.\\
By Lemma 2.1(1) and Lemma 2.1(6), we have
\begin{equation*}
\begin{split}
n'&\leq A_1+A_2+\cdots+A_k-A_k\\
& = A_1+A_2+\cdots+A_{k-1}\\
& < A_{k+1}\\
n'& \geq A_1+A_2+\cdots+A_k - A_{k,0}\\
& > A_{k,0}\\
\end{split}
\end{equation*}

Hence $P_{n'}$ can be determined by the propositions above.
Moreover, by Proposition 4.2 $U_n$ and $D_n$ are dual to $D_{n'}$
and $U_{n'}$, respectively. So they can also be determined by
recursions.

\begin{proof}[Proof of Proposition 4.3]
Assume $S_n \neq \emptyset$. (Otherwise the proposition is trivially
true.) Let $\hat{1}$ be the maximal element in $P_n$. Then
$l(\hat{1}) = k$ (Proposition 3.1). So $\hat{1}-\tau_k \in
S_{n-A_k}$, which implies $S_{n-A_k} \neq \emptyset$. Note that
\[n-A_k > A_{k,2}-A_k = \sum_{\substack{1\leq i < k-N+1 \\ N\nmid
i}}A_{k-N+1-i}.\] By Corollary 3.2, we have
$n-A_k \geq A_{k-N+1}$.\\
For any $a \in S_n$,
\begin{equation*}
\begin{split}
s(a)=n&\geq A_k+A_{k-N+1}\\
& > (A_{k-1}+\cdots+A_{k-N})+(A_1+A_2+\cdots+A_{k-N-1})
(\textrm{Lemma }
2.1(1))\\
& = A_1+A_2+\cdots+A_{k-1}\\
\end{split}
\end{equation*}
Hence $l(a) = k$ for all $a \in S_n$. In particular, it is valid to
do the subtraction $a - \tau_k$.\\
Therefore, the map
\begin{equation*}
\begin{split}
\phi : P_n & \rightarrow P_{n-A_k}\\
a & \mapsto a - \tau_k\\
\end{split}
\end{equation*}
gives an isomorphism from $P_n$ to $P_{n-A_k}$ and $P_n =
P_{n-A_k}+\tau_k$.
\end{proof}

\begin{proof}[Proof of Proposition 4.4]
Assume $S_n \neq \emptyset$. In this case, $S_n = T_n \cup
R_n$.\\Applying a similar argument, we obtain $n \geq A_k+A_{k-N}$.
Furthermore, by Proposition 3.1,
$l(\hat{1}-\tau_k) = k-N$.\\
Define
\begin{equation*}
\begin{split}
\phi_1: U_n& \rightarrow P_{n-A_k}\\
a& \mapsto a-\tau_k\\
\end{split}
\end{equation*}
Then $\phi_1$ gives an isomorphism from $U_n$ to $P_{n-A_k}$ and
$U_n = P_{n-A_k}+\tau_k$.\\
On the other hand, for any $a=(a_1,a_2, \cdots, a_{k-1}) \in D_n$,
we claim that $a_{k-N} = a_{k-N+1}= \cdots = a_{k-1}$.\\
Assume the contrary. Then
\begin{equation*}
\begin{split}
s(a) & \leq A_1 +\cdots+A_{k-N-1}+A_{k-N+1}+\cdots+A_{k-1}\\
& < (A_{k-N}+A_{k-N-1}) +A_{k-N+1}+\cdots+A_{k-1}\\
& = A_{k-N-1}+A_k\\
& < A_{k-N}+A_k\\
&\leq n\\
\end{split}
\end{equation*}
a contradiction.\\
Hence, the subtraction $a-\eta_k$ is valid.\\
Note that $l(a-\eta_k)\leq k-N-1 = l(\hat{1}-\tau_k)-1$. Thus
$a-\eta_k \in D_{n-A_k}$ for all $a\in D_n$. So the map
\begin{equation*}
\begin{split}
\phi_2: D_n & \rightarrow D_{n-A_k}\\
a & \mapsto a-\eta_k\\
\end{split}
\end{equation*}
gives an isomorphism from $D_n$ to $D_{n-A_k}$ and $D_n =
D_{n-A_k}+\eta_k$.\\
Therefore, \[P_n = (P_{n-A_k}+\tau_k) \hat{+}(D_{n-A_k}+\eta_k)\]
via the natural map $\psi: D_{n-A_k}+\tau_k \rightarrow
D_{n-A_k}+\eta_k$: $a+\tau_k \mapsto a+\eta_k$.
\end{proof}

\begin{proof}[Proof of Proposition 4.5]
The proof is almost the
same. Assume $S_n \neq \emptyset$.\\
In this case, we have $n \geq A_k + A_{k-N-1}$ and
$l(\hat{1}-\tau_k)=k-N-1$.\\
For $a\in D_n$, we still have $a_{k-N}=a_{k-N+1}=\cdots=A_{k-1}=1$.
So the map
\begin{equation*}
\begin{split}
\phi_2: D_n & \rightarrow P_{n-A_k}\\
a & \mapsto a-\eta_k\\
\end{split}
\end{equation*}
gives an isomorphism from $D_n$ to $P_{n-A_k}$. This completes the
proof.
\end{proof}

\section{Applications}

\subsection{$P_n$ are modular lattices}

\begin{thm}
$P_n$, $U_n$, $D_n$ are modular lattices.
\end{thm}

\begin{proof}
We prove the statement by induction on $n$.\\
\textbf{Base case:} When $n < A_{N+1}$, $P_n$ is either $\emptyset$
or a single-element set. The statement is trivially true.\\
\textbf{Inductive step:} Consider $A_k \leq n < A_{k+1}$ ($k\geq
N+1$). It is clear that the dual poset of a modular lattice is also
a modular lattice. So, by symmetry, it suffices to consider
$A_{k,0}< n <A_{k+1}$.
\begin{enumerate}
\item If $A_{k,2}<n<A_{k+1}$, then $P_n$, $U_n$, $D_n$ are isomorphic
to $P_{n-A_k}$, $U_{n-A_k}$ and $D_{n-A_k}$, respectively. So they
are modular lattices by the induction hypothesis.
\item If $A_{k,1}<n\leq A_{k,2}$, then $U_n$, $D_n$ are isomorphic
to $P_{n-A_k}$ and $D_{n-A_k}$, respectively. So they are lattices
by induction. To show that $P_n$ is a lattice, we need to show that
for any
$x,y\in P_n$, $x\vee y$ and $x\wedge y$ exist.\\
Indeed, if $x,y\in U_n$ or $x,y\in D_n$, then $x\vee y$ and $x\wedge
y$ exist by the induction hypothesis.\\
If $x\in U_n$ and $y\in D_n$, it is easy to check that $x\vee y$ is
simply $x\vee\psi^{-1}(y)$ in $P_{n-A_k}$. Similarly, $x\wedge y$ is
simply $\psi(x)\wedge y$ in $D_{n-A_k}$. Hence, $P_n$ is a
lattice.\\
To show modularity, we need to show that $x,y$ both cover
$x\wedge y$ if and only if $x$,$y$ are both covered by $x\vee y$.\\
Indeed, if $x,y \in U_n$ and $x,y \in D_n$, the statement follows by
the induction hypothesis.\\
Now suppose $x \in U_n$ and $y \in D_n$. Without loss of generality,
assume that $x\vee y$ covers both $x$,
$y$.\\
Obviously, $x\vee y \in U_n$. So $x\vee y = \psi^{-1}(y)$ and
$\psi^{-1}(y)$ covers $x$.\\
Therefore, $\psi(x)$ is covered by both $x$ and $y$, as desired.
\item The case $A_{k,0}<n\leq A_{k,1}$ is similar to case(2).

\end{enumerate}

\end{proof}

\subsection{Quasifibonacci Sequence of Even Level}
In this section we will prove Theorem 1.1. As mentioned in the
introduction section, $h_m$ is the difference of the number of
partitions into an even number of $A_k$'s and the number of
partitions into an odd number of $A_k$'s.\\
To distinguish these two kinds of partitions, we define the
\emph{sign function} $\sigma: \{0,1\}^{\omega} \rightarrow \pm1$ by
setting $\sigma(a) = 1$ if $a$ contains an even number of 1's and
$\sigma(a) = -1$ otherwise. In general, for any finite subset $P$ of
$\{0,1\}^{\omega}$ with natural partial ordering, define
\[\sigma(P) := \sum_{a\in P} \sigma(a)\] It is clear that $h_n = \sigma(P_n)$. We also define $f_n =\sigma(U_n)$ and
$g_n = \sigma(D_n)$. Then, obviously, $h_n = f_n+g_n$.

The following lemma will be useful in the proof below.
\begin{lem}
\begin{itemize}
\item Let $a,b\in \{0,1\}^{\omega}$. Then $\sigma(a\pm b) = \sigma(a)\sigma(b)$
\item Let $a \in \{0,1\}^{\omega}$ and $P$ a finite subset of $\{0,1\}^{\omega}$. Then
$\sigma(P+a)=\sigma(P)\sigma(a)$
\item Let $P$, $Q$ be finite subsets of $\{0,1\}^{\omega}$. Then $\sigma(P\hat{+}Q) =
\sigma(P)+\sigma(Q)$.
\end{itemize}
\end{lem}

\begin{proof}
Straightforward.
\end{proof}

To visualize the relation of odd and even partitions, we color the
digraph with two colors. In the corresponding Hasse diagram, a
vertex $a \in P_n$ is colored blue if $\sigma(a) = 1$ and colored
red if $\sigma(a) = -1$. The figure below shows the first 24 colored
Hasse diagrams for $P_n(\{F_k\})$:
\begin{figure}[h]

  \begin{center}
  \psfrag{p1}{$P_1$}
  \psfrag{p2}{$P_2$}
  \psfrag{p3}{$P_3$}
  \psfrag{p4}{$P_4$}
  \psfrag{p5}{$P_5$}
  \psfrag{p6}{$P_6$}
  \psfrag{p7}{$P_7$}
  \psfrag{p8}{$P_8$}
  \psfrag{p9}{$P_9$}
  \psfrag{p10}{$P_{10}$}
  \psfrag{p11}{$P_{11}$}
  \psfrag{p12}{$P_{12}$}
  \psfrag{p13}{$P_{13}$}
  \psfrag{p14}{$P_{14}$}
  \psfrag{p15}{$P_{15}$}
  \psfrag{p16}{$P_{16}$}
  \psfrag{p17}{$P_{17}$}
  \psfrag{p18}{$P_{18}$}
  \psfrag{p19}{$P_{19}$}
  \psfrag{p20}{$P_{20}$}
  \psfrag{p21}{$P_{21}$}
  \psfrag{p22}{$P_{22}$}
  \psfrag{p23}{$P_{23}$}
  \psfrag{p24}{$P_{24}$}
  \includegraphics[width=12cm]{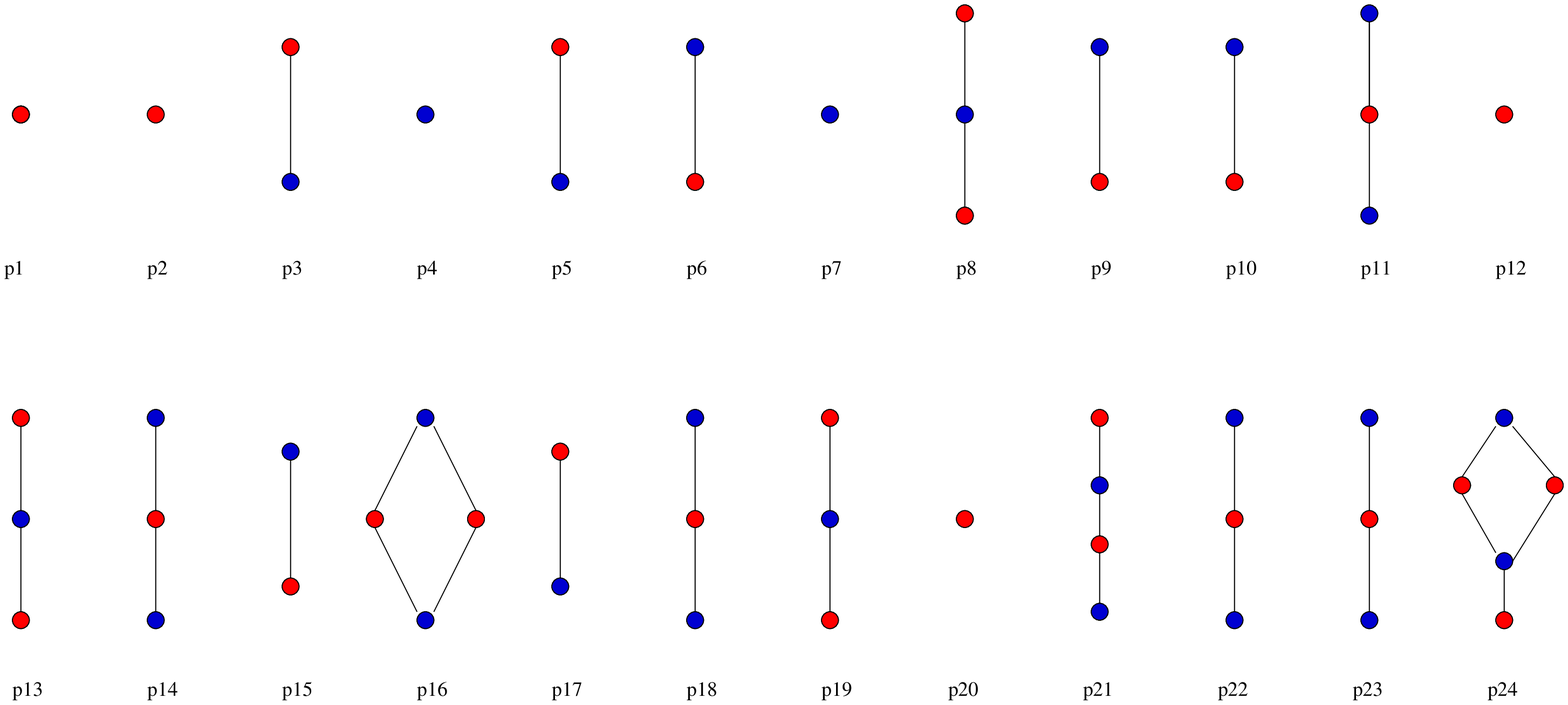}\\

  \end{center}
\end{figure}

Suppose $a$ covers $b$ in $P_n$. By the definition, $b$ has $N-1$
more 1's than $a$. Hence adjacent vertices in the corresponding
digraph have different colors if $N$ is even.

Now we prove Theorem 1.1 by proving the following stronger result:

\begin{prop}
For any $n \in \mathbb{Z}^+$, $f_n$, $g_n$, $h_n \in \{-1, 0, 1\}$.

\end{prop}

\begin{proof}
We perform induction on $n$.\\
\textbf{Base case: }When $n < A_{N+1}$, $P_n$ is either $\emptyset$
or a single-element set. The statement is true.\\
\textbf{Inductive step:} Consider $A_k \leq n < A_{k+1}$ $(k\geq
N+1)$.
\begin{enumerate}
\item If $A_{k,2}<n<A_{k+1}$, then $U_n = U_{n-A_k} + \tau_k$, $D_n = D_{n-A_k} + \tau_k$, $P_n = P_{n-A_k}+\tau_k$. Hence we
have
\[f_n = \sigma(\tau_k) f_{n-A_k} = -f_{n-A_k} \in \{-1,0,1\}\]
Similarly, we have $g_n = -g_{n-A_k} \in \{-1,0,1\}$ and $h_n =
-h-{n-A_k} \in \{-1,0,1\}$.
\item If $A_{k,1} < n \leq A_{k,2}$, then $U_n = P_{n-A_k}+\tau_k$
and $D_n = P_{n-A_k}+\eta_k$ and $P_n =
(P_{n-A_k}+\tau_k)\hat{+}(D_{n-A_k}+\eta_k)$. Hence we have
\[f_n = \sigma(\tau_k) h_{n-A_k} = -h_{n-A_k} \in \{-1,0,1\}\]
\[g_n = \sigma(\eta_k) g_{n-A_k} = g_{n-A_k} \in \{-1,0,1\}\]
and
\begin{equation*}
\begin{split}
h_n& = \sigma(P_{n-A_k}+\tau_k) + \sigma(D_{n-A_k}+\eta_k)\\
& = \sigma(P_{n-A_k})\sigma(\tau_k) +
\sigma(D_{n-A_k})\sigma(\eta_k)\\
& = -h_{n-A_k} + g_{n-A_k}\\
& = -f_{n-A_k} \in \{-1, 0, 1\}\\
\end{split}
\end{equation*}
\item If $A_{k,0} < n \leq A_{k,1}$, then $U_n = P_{n-A_k}+\tau_k$,
$D_n=P_{n-A_k}+\eta_k$ and $P_n = (U_{n-A_k}+\tau_k)
\hat{+}(D_{n-A_k}+\eta_k)$. Hence we have
\[f_n = \sigma(\tau_k) \sigma(P_{n-A_k}) = -h_{n-A_k} \in \{-1,0,1\}\]
\[g_n = \sigma(\eta_k) \sigma(P_{n-A_k}) = h_{n-A_k} \in \{-1,0,1\}\]
and
\begin{equation*}
\begin{split}
h_n & =\sigma(P_{n-A_k}+\tau_k)+\sigma(P_{n-A_k}+\eta_k)\\
& = \sigma(P_{n-A_k})\sigma(\tau_k) +
\sigma(P_{n-A_k})\sigma(\eta_k)\\
& = -h_{n-A_k}+h_{n-A_k}\\
& = 0\\
\end{split}
\end{equation*}
\item If $A_k\leq n \leq A_{k,0}$, then $U_n$, $D_n$, $P_n$ are
dual to $D_{n'}$, $U_{n'}$ and $P_{n'}$, respectively, where $n' =
A_1+A_2+\cdots+A_k - n$. Hence we have
\[f_n = (-1)^k \sigma(D_{n'}) =(-1)^k g_{n'}\in \{-1,0,1\}^{\omega} \]
\[g_n = (-1)^k \sigma(U_{n'}) =(-1)^k f_{n'}\in \{-1,0,1\}^{\omega} \]
and \[h_n = (-1)^k \sigma(P_{n'}) =(-1)^k h_{n'}\in
\{-1,0,1\}^{\omega}\] This completes the proof.

\end{enumerate}
\end{proof}

To write down the recursion explicitly, we have $(f_n, g_n, h_n)=$
\[\left\{ \begin{array}{llll} (-1)^{k+1} g_{\gamma_k-n-A_k} &
(-1)^{k+1} f_{\gamma_k-n-A_k} & (-1)^{k+1} h_{\gamma_k-n-A_k} &
\mbox{if $A_k \leq n < \gamma_k -
A_{k,2}$};\\
(-1)^k g_{\gamma_k-n-A_k} & (-1)^{k+1} h_{\gamma_k-n-A_k} &
(-1)^{k+1} f_{\gamma_k-n-A_k} & \mbox{if $ \gamma_k - A_{k,2} \leq n
< \gamma_k -
A_{k,1}$};\\
(-1)^k h_{\gamma_k-n-A_k} &(-1)^{k+1} h_{\gamma_k-n-A_k} & 0 &
\mbox{if $ \gamma_k - A_{k,1} \leq n
\leq A_{k,0}$};\\
-h_{n-A_k} &  h_{n-A_k} & 0 & \mbox{if $A_{k,0} < n
\leq A_{k,1}$};\\
-h_{n-A_k} &  g_{n-A_k} &  -f_{n-A_k} & \mbox{if $A_{k,1} < n
\leq A_{k,2}$};\\
-f_{n-A_k} & -g_{n-A_k} &  -h_{n-A_k} & \mbox{if $A_{k,2} < n
<A_{k+1}$}. \end{array} \right.\] where $\gamma_k =
A_1+A_2+\cdots+A_k$.

\emph{Comment:} If we take $A_k = F_k$, then the recursion above is
precisely the one in [1].

\subsection{Quasifibonacci Sequence of Odd Level}
In this case, $N$ is an odd number. Hence for each pair of adjacent
vertices $(a, b)$ in $G_n$, the parity of the number of 1's in $a$
and $b$ must be the same. In other words, every adjacent pair of
vertices have same color. Since $P_n$ are lattices, $G_n$ are
connected graphs. Therefore all vertices in $P_n$ have same color.
It is easily seen that $h_n$ is not bounded in this case. Instead we
have the following estimate.

\begin{prop}
Let $k \geq N$ be an integer. For any $A_k\leq n <A_{k+1}$, $|h_n|
\leq 2^{k-N}$.

\end{prop}

\begin{proof}
\textbf{Base case} When $n < A_{N+1}$, $P_n$ is either $\emptyset$
or a single-element set. So $h_n \in \{0,-1\}$. The statement is
true.\\
\textbf{Inductive step} Consider $A_k \leq n <A_{k+1}$ $(k \geq
N+1)$. With a similar argument, we can derive the following
recursion:
\[h_n = \left\{ \begin{array}{ll} (-1)^{k+1} h_{\gamma_k-n-A_k} &
\mbox{if $A_k \leq n < \gamma_k -
A_{k,2}$};\\
(-1)^{k+1} (g_{\gamma_k-n-A_k}+h_{\gamma_k-n-A_k}) & \mbox{if $
\gamma_k - A_{k,2} \leq n < \gamma_k -
A_{k,1}$};\\
(-1)^{k+1}2h_{\gamma_k-n-A_k} & \mbox{if $ \gamma_k - A_{k,1} \leq n
\leq A_{k,0}$};\\
-2h_{n-A_k} & \mbox{if $A_{k,0} < n
\leq A_{k,1}$};\\
-g_{n-A_k}-h_{n-A_k} & \mbox{if $A_{k,1} < n
\leq A_{k,2}$};\\
-h_{n-A_k} & \mbox{if $A_{k,2} < n <A_{k+1}$}. \end{array} \right.\]
where $\gamma_k = A_1+A_2+\cdots+A_k$.\\
Note that $|g_m| \leq |h_m|$ $(\forall m \in \mathbb{Z}^+)$. Hence,
in any of the six cases, $|h_n| \leq 2\cdot 2^{k-N-1} = 2^{k-N}$.
This completes the proof.
\end{proof}
\emph{Comment:} This upper bound is the best possible because there
exists $A_k \leq n <A_{k+1}$ satisfying $|h_n| = 2^{k-N}$ for all $k
\geq N$. However, it is possible to improve the result by spliting
the intervals into more pieces and refining the estimate.

\section*{Acknowledgements}
The author would like to thank professor Richard Stanley for
valuable suggestions on this paper.

\bigskip

\footnotesize\textsc{Department of Mathematics, Massachusetts
Institute of Technology, Cambridge, MA 02139}

\end{document}